\documentclass[12pt]{article}

\usepackage{amsfonts,amssymb,amsmath}
\numberwithin{equation}{section}

\newcommand{\Z}{\mathbb{Z}}
\newcommand{\C}{\mathbb{C}}

\newcommand{\la}{\lambda}

\newcommand{\bK}{\mathcal{K}}

\newcommand{\fM}{\mathfrak{M}}

\newcommand{\vph}{\varphi}
\newcommand{\tV}{V^*}

\newcommand{\z}{\zeta}
\newcommand{\w}{\eta}
\newcommand{\V}{\mathfrak{V}}

\newtheorem{Theorem}{Theorem}

\newtheorem{Remark}{Remark}
\newtheorem{Example}{Example}

\begin{document}

\title{A Fredholm determinant formula for Toeplitz determinants}
\author{Alexei Borodin and Andrei Okounkov}
\date{}
%\thanks{We would like to thank MSRI and the organizers of the
%Random}

\maketitle

\section{Introduction}

The purpose of this note is to explain how the results of \cite{O2}
apply to a question raised by A.~Its and, independently, P.~Deift
during the MSRI workshop on Random Matrices in June 1999.
The question was whether there exists a general formula expressing a
Toeplitz determinant
$$
D_n(\vph)=\det\big(\vph_{i-j}\big)_{1\le i,j\le n}\,,
\quad \vph(\z)=\sum_{k\in\Z} \vph_k \, \z^k\,,
$$
as the Fredholm determinant of an
operator $1-\bK$ acting on $\ell_2(\{n,n+1,\dots\})$,
where the kernel $\bK=\bK(\vph)$ admits an integral representation
in terms of $\vph$. The answer is affirmative and the construction
of the kernel is explained below.

We give two versions of the result: an algebraic one, which holds
in the suitable algebra of formal power series, and an analytic one.
In order to minimize the amount of analysis, we make a rather restrictive 
analyticity assumption on the function $\vph$. One should be able to
relax this assumption, see Remark \ref{weak}. Our proof is a
direct application of two results due to I.~Gessel \cite{G} 
and one of the authors \cite{O2}, respectively.

We also consider 3 examples in which
  the kernel $\bK$ can be expressed in 
classical special functions. First two of them have been worked out
in \cite{BDJ1,TW,BOO,J}.

\section{Statement of the result}

We assume that $\vph$ is an exponential of a Laurent series
$\vph(\z)=\exp{V(\z)}$ and we assume that $\vph$ is normalized in
such a way that $V(\z)$
has no constant term. Let us separate
positive and negative
powers in this series
$$
V(\z)=V^+(\z)+V^-(\z)\,,
$$
where $V^\pm(\z)=\sum_{k=1}^\infty v_k^{\pm}\z^{\pm k}$.

First, consider the algebraic situation.  
Introduce a $\Z_+$--grading in the polynomial
algebra $\C[v^{\pm}_k]_{k\ge 1}$ by setting
$$
\deg v^{\pm}_k = k\,.
$$
For an arbitrary $f\in \C[v^{\pm}_k]_{k\ge 1}$ we denote by
$\deg f$ the \emph{minimal} degree of all terms in $f$. 
Let $\V$ be the completion of $\C[v^{\pm}_k]_{k\ge 1}$ with
respect to this grading. Recall that $f_n\to f$ in $\V$ by definition
means that
$$
\deg (f-f_n) \to \infty \,, \quad n\to\infty\,.
$$

Because $V$ has no constant term
it follows that $\vph(\z)=\exp{V(\z)}$
is a well-defined Laurent series in $\z$ with coefficients in $\V$ and,
moreover, 
$$
\deg\, \big[\z^k\big] \,\vph(\z) \ge |k|  \,,
$$
where $\big[\z^k\big] \vph(\z)$ stands for the coefficient 
of $\z^k$ in $\vph$. 

By definition, set
$$
\tV(\z)= V^-(-\z)- V^+(-\z) \,,
$$
and define the kernel $\bK$ by the following generating function
\begin{equation}\label{e0}
\sum_{i,j\in\Z} \z^i \w^{-j} \, \bK(i,j) = \frac{
\exp
\left(\tV(\z)-
\tV(\w)\right)
}
{\z/\w-1}  \,, \quad |\z|>|\w| \,.
\end{equation}
Here the condition $|\z|>|\w|$ means that $\frac1{\z/\w-1}$ is 
to be expanded in powers of $\frac{\w}{\z}$. In other words,
$$
\bK(i,j)=\sum_{l\ge 1} \big[\z^{i+l}\, \w^{-j-l}\big] \, \exp
\left(\tV(\z)-
\tV(\w)\right) \,,
$$
where $\big[\z^{i+l}\, \w^{-j-l}\big]$ stands for the coefficient
of $\z^{i+l}\, \w^{-j-l}$. Since 
$$
\deg \big[\z^{k}\big] \,  \exp
\left(\tV(\z)\right) \ge |k|
$$
we conclude that $\bK(i,j)$, $i,j\ge 0$, is a well defined
element of $\V$ and
$$
\deg \, \bK(i,j) \ge i+j+2 \,.
$$

In the analytic situation, we will assume
that the Laurent series $V(\z)$ converges in some
nonempty annular neighborhood of the unit circle
$$
r^{-1}<|\z|<r \,,\quad \z\in\C\,,\quad r>1 \,.
$$
In this case the kernel $\bK$ can be expressed
by the following contour integral
\begin{equation}\label{e1}
\bK(i,j)=
\frac1{(2\pi\sqrt{-1})^2} \iint_{|\z|>|\w|}
\frac{
\exp
\left(\tV(\z)-
\tV(\w)\right)
}
{\z-\w}\,
\frac{d\z}{\z^{i+1}}\, \frac{d\w}{\w^{-j}} \,.
\end{equation}
Here the integral is taken, for example,
over the torus $|\z|=|\w|^{-1}=\rho$, $1<\rho<r$, in $\C^2$.
It is clear that
$$
\sum_{i,j\ge 0} |\bK(i,j)| \le \textup{const}\, 
\sum_{i,j\ge 0} \frac1{\rho^{i+j+1}} < \infty 
$$
and therefore $\bK$ defines a trace class operator on 
$\ell_2(\{0,1,\dots,\})$.

Finally, set
$$
Z= \exp\left(\sum_{k=1}^\infty k\, v^+_k \, v^-_{k} \right)\,.
$$
This is easily seen to be well defined in both algebraic and
analytic case. Our main result is the following 

\begin{Theorem}\label{main} We have the following identity
\begin{equation}\label{algf}
D_n\left(\vph\right) =Z\,  \det (1-\bK)_{\ell_2(\{n,n+1,\dots\})}\,,
\end{equation}
where $D_n\left(\vph\right)$ is the
$n\times n$ Toeplitz determinant with symbol $\vph(\z)
=\exp{V(\z)}$ and the Fredholm determinant in the right-hand side
is defined by 
\begin{equation}\label{fred}
\det (1-\bK)_{\ell_2(\{n,n+1,\dots\})}=\sum_{m=0}^\infty
{(-1)^m}\sum_{n\le l_1<\dots<l_m}^\infty\det\,\big[\bK(l_i,l_j)
\big]_{i,j=1}^m.
\end{equation}
The convergence in \eqref{fred} is the usual convergence in 
the case when $V(\z)$ converges in some neighborhood of the
unit cirle in $\C$ or else 
the convergence in the algebra $\V$ of formal power series. 
\end{Theorem}

\begin{Remark}\label{szego}\rm The $n\to\infty$ limit
corresponds to the situation of the strong Szeg\"o limit theorem \cite{sz}
$$
D_n(\vph)\to Z, \qquad n\to\infty\,.
$$
In the algebraic situation, this convergence means that
for any fixed $k$ the terms of degree $\le k$ in $D_n$ 
will coincide with those of $Z$ for large enough $n$. 
\end{Remark}

\begin{Remark}\label{integr} \rm
Replacing $\z$ and $\w$ in the right--hand side of
\eqref{e1} by $t \z$ and $t \w$, respectively,
and applying
$\left.\frac{d}{dt}\right|_{t=1}$ yields
\begin{multline}\label{e2}
(i-j)\, \bK(i,j)\\=
 \frac1{(2\pi\sqrt{-1})^2}
\iint_{|\z|>|\w|}
\frac{\z\dfrac{d}{d\z} \tV(\z)- \w\dfrac{d}{d\w} \tV(\w)}{\z-\w} \,
e^{\tV(\z)-
\tV(\w)}\, \frac{d\z}{\z^{i+1}} \, \frac{d\w}{\w^{-j}} \,.
\end{multline}

The advantage of the formula \eqref{e2} over the formula \eqref{e1}
is that often the division in \eqref{e2} can be carried out explicitly and
this leads to an expression  of the type
$$
 (i-j)\, \bK(i,j) = \sum_{k=1}^m f_k(i) \, g_k(j) \,,
$$
for some $m$ and some functions  $f_1,\dots, f_m$ and $g_1,\dots,g_m$.
This means that the kernel $\bK$ is \emph{integrable}
in the sense of \cite{IIKS}, see also \cite{D1}.
\footnote{Note, however, that the notion of an integrable operator
was introduced and extensively used for kernels with \emph{continuous}
arguments while in our case the set of arguments is \emph{discrete}.}

The main disadvantage of the formula \eqref{e2} is that it does not
allow to compute the diagonal entries $\bK(i,i)$.
\end{Remark}

\begin{Remark}\label{weak}\rm 
The analyticity assumption on $V(z)$ can, presumably, be 
seriously weakened. One obvious remark is that if we have
a sequence of convergent Laurent series $V_n(\z)$ such
that the coefficients of $e^{V_n(z)}$ converge as $n\to\infty$ and the 
associated kernels $\bK_n$ converge in the trace norm
in $\ell_2(\{0,1,\dots\})$, then Theorem \ref{main} is
satisfied for the limit series $V(z)$ which may not
represent an analytic function. 
\end{Remark}  

\begin{Remark}\label{ind}\rm 
Finally notice that the left-hand side in \eqref{algf}
obviously depends only on $\vph_{-n+1},\dots,\vph_{n-1}$,
which is not so clear from the right-hand side. 
\end{Remark}

\section{Proof}

The proof of Theorem \ref{main} in the analytic and 
algebraic setup is identical, with the understanding that
the meaning of convergence of infinite series 
in these two cases is very different. 
Theorem \ref{main} will be established by applying
the results of \cite{O2} to
a formula, due to Gessel \cite{G}, expressing $D_n(\vph)$
as a series in Schur functions. Although neither side of the
equality in the theorem involves Schur functions, they will be
needed in our proof.

Recall \cite{M}, I.3,
that Schur functions $s_\la$  are certain very distinguished
symmetric polynomials in auxiliary variables
$x=(x_1,x_2,\dots)$. For our purposes and notational
compatibility with \cite{O2},  it is more
convenient to regard $s_\la$ as polynomials in the power-sum symmetric
functions
$$
t_k = k^{-1} \sum x_i^k  \,, \quad k=1,2,\dots \,,
$$
which are free commutative generators of the algebra of symmetric functions.
Since we will never use the original variables $x$, we will
abuse the notation and write $s_\la(t)$ for the polynomial in the $t_k$'s
which gives the corresponding Schur function.
We do not need the explicit form of these
polynomials but it is worth mentioning that their coefficients
are the characters of symmetric groups, see \cite{M}, Section I.7.

We  now recall a result of Gessel \cite{G} or, more precisely,
its dual version \cite{TW}.
Consider the following sum:
$$
\sum_{\la,\, \la_1\le n} s_\la(t^+)\, s_\la(t^-)\,,
$$
where
$$
t^\pm=(t^\pm_1,t^\pm_2,\dots)
$$
are two sets of variables and
the sum is taken over all
partitions $\la$ with first part $\la_1$ less or equal to $n$.

Gessel's theorem
asserts that
\footnote{
To see the connection with Section 4 of \cite{TW} note that
if $t_k = k^{-1} \sum x_i^k$ then $\sum_{k\ge 1} t_k \z^k =
\prod_i (1-x_i \z)^{-1}$.
}
\begin{equation}\label{gess}
 \sum_{\la, \, \la_1\le n} s_\la(t^+)\, s_\la(t^-) =
D_n\left(\exp\left[-T^+(-\z)-T^-(-\z)\right]\right)
\end{equation}
where $T^\pm(z)=\sum_{k=1}^\infty t^\pm_k \, \z^{\pm k}$. Let us
identify the parameters $v^\pm_k$ and $t^\pm_k$ by setting
\begin{equation}\label{e6}
V(\z)= -T^+(-\z)-T^-(-\z) \,.
\end{equation}
That is, $v^\pm_k = (-1)^{k+1}\, t^\pm_k$ \,.

Following \cite{O2}, consider the following 
\emph{Schur measure} on partitions 
$$
\fM(\la)=\frac1Z\, s_\la(t^+)\, s_\la(t^-) \,
$$
where $Z$ is the sum in the Cauchy identity for the Schur
functions
$$
Z= \sum_\la s_\la(t^+)\, s_\la(t^-) =
\exp\left(\sum_k k \, t^+_k\, t^-_k\right) =
\exp\left(\sum_{k=1}^\infty k\, v^+_k \, v^-_{k} \right)\,.
$$
Note that although we call $\fM$ a measure it does not have
to be positive or even real for our purposes. It is, however,
positive if $t^-_k = \overline{t^+_k}$, $k=1,2,\dots$, which
is the case if $V(z)$ is real on the unit circle.

Given a partition $\la$, set $S(\la)=\{\la_i-i\}\subset\Z$.
The Gessel theorem
can be rewritten as follows
$$
D_n\left(\vph\right) =
Z \, \sum_{\la, \, \la_1\le n}\fM(\la)\notag
=Z \, \sum_{\la, \, S(\la)\cap\{n,n+1,\dots\}=
\varnothing}\fM(\la)\label{e4}\,.
$$
It was proved in \cite{O2}, Theorems 1 and 2,
that for any finite subset
$X\subset\Z$
\begin{equation}\label{e5}
\sum_{\la, \, X\subset S(\la)}\fM(\la) =
\det\Big[\bK(x_i,x_j)\Big]_{x_i,x_j\in X} \,.
\end{equation}
Here the kernel $\bK$ is given by the
generating function
\footnote{Note that our kernel $\bK$ differs from the
kernel $K$ used in \cite{O2} by a shift of variables
$\bK(x,y)=K(x+\frac12,y+\frac12)$.

Also note that we identify the variables $t^\pm$ with the variables
$t$ and $t'$ used in \cite{O2} as follows: $t=t^-$, $t'=t^+$.}
$$
\sum_{i,j\in\Z} \z^i \w^{-j} \, \bK(i,j)=
\frac{\exp\big(T^+(\z)-T^-(\z)-T^+(\w)+T^-(\w)\big)}{\z/\w-1}\,,
\quad |\w|<|\z|\,,
$$
which is immediately seen to be identical to \eqref{e0}
after the identification \eqref{e6}.

By the usual inclusion--exclusion principle we have
\begin{multline*}
\fM\left(\left\{\la, \, S(\la)\cap\{n,n+1,\dots\}=
\varnothing\right\}\right)=\\
\sum_{m=0}^{\infty} (-1)^m \sum_{n\le l_1<\dots<l_m}^\infty
\fM\left(\left\{\la, \, \{l_1,\dots,l_m\}\subset S(\la)
\right\}\right)
\end{multline*}
This formula together with \eqref{e5} concludes the proof of \eqref{algf}.

\section{Examples}
\begin{Example}\label{example1} \rm We start with the simplest nontrivial example when
$v_k^\pm=0$ for all $k>1$. Then $\vph(\z)=\exp(v_1^+\z+v_1^-\z^{-1})$.
Since Toeplitz determinants $D_n(\vph)$ do not change under the transformation
$\vph(\z)\to \vph(a\z)$, $a\in \C\setminus\{0\}$, all that matters in our case
is the product $v_1^+v_1^-$.  Thus, we may assume that
$v_1^+=v^-_1=\theta\in \C\setminus\{0\}$. Note that the asymptotics of
Toeplitz determinants with symbol $\vph(\z)=\exp(\theta(\z+\z^{-1}))$ played
a central role in \cite{BDJ1}.
We have
$$
V(\z)=\theta(\z+\z^{-1}),\quad \tV(\z)=\theta(\z-\z^{-1}).
$$
To compute the kernel $\bK$ we will employ \eqref{e2}. We get
$$
\frac{\z\dfrac{d}{d\z} \tV(\z)- \w\dfrac{d}{d\w} \tV(\w)}{\z-\w} =\theta-
\frac{\theta}{\z\w},
$$

\begin{multline}\label{bess}
(i-j)\,\bK(i,j)\\ =\frac\theta{(2\pi\sqrt{-1})^2}
\iint\limits_{|\z|=c_1,\,|\w|=c_2}
\Bigl(
\frac{e^{\theta(\z-\z^{-1})}}{z^{i+1}}\,\frac
{e^{\theta(\w^{-1}-\w)}}{\w^{-j}}-
\frac{e^{\theta(\z-\z^{-1})}}{z^{i+2}}\,\frac
{e^{\theta(\w^{-1}-\w)}}{\w^{-j+1}}\Bigr)\,d\z d\w
\end{multline}

Since
$$
\exp(\theta(\xi-\xi^{-1}))=\sum_{k=-\infty}^\infty \xi^kJ_k(2\theta),
$$
where $J_{\nu}(x)$ is the Bessel function, see \cite{HTF}, 7.2.4,
we can rewrite \eqref{bess} as follows
\begin{equation}\label{planch}
\bK(i,j)=\theta\,\frac{J_i(2\theta)J_{j+1}(2\theta)-J_{i+1}(2\theta)
J_{j}(2\theta)}{i-j}.
\end{equation}

It can be shown, for example, by expanding $(\z-\w)^{-1}$ in \eqref{e1}
as the sum of a geometric progression, that the diagonal
entries $\bK(i,i)$ are given by the limit of the right--hand side of
\eqref{planch} as $j\to i$. The l'Hopital rule yields
$$
\bK(i,i)=\theta\,\left(\frac d{di}J_i(2\theta)\,J_{i+1}(2\theta)-
\frac d{di}J_{i+1}(2\theta)\,J_i(2\theta)\right).
$$

The kernel $\bK$ has been obtained in \cite{BOO} and,
independently, in \cite{J}
in connection with asymptotics of the Plancherel measures for
symmetric groups.

The constant $Z$ is this example equals $\exp(\theta^2)$, and
Theorem \ref{main} reads
\begin{equation}\label{ex1}
D_n(\exp(\theta(\z+\z^{-1}))=\exp(\theta^2)
\,\det (1-\bK)_{\ell_2(\{n,n+1,\dots\})}\,
\end{equation}
where $\bK$ is as above.

This formula also follows from the results
of \cite{G,BOO,J}.
\end{Example}

\begin{Example}\label{example2} \rm 
In this example we will consider Toeplitz determinants with symbols
of the form $\vph(\z)=(1+\theta_1\z)^\kappa\exp(\theta_2/\z)$. Similarly to
the previous example, since dilations do not change Toeplitz determinants,
we may assume that $\theta_1=\theta_2=\theta$. Then
$$
V(\z)=\ln\vph(\z)=\kappa\ln(1+\theta\z)+\theta\z^{-1}.
$$
Note that in order to satisfy
the conditions of Theorem \ref{main}
we have to assume that $|\theta|<1$.
We have
$$
\tV(\z)=-\kappa\ln(1-\theta\z)-\theta\z^{-1},\quad
\z\frac d{d\z}\tV(\z)=\frac{\kappa\theta\z}{1-\theta\z}+\theta\z^{-1},
$$
$$
\frac{\z\dfrac{d}{d\z} \tV(\z)- \w\dfrac{d}{d\w} \tV(\w)}{\z-\w}=
\frac{\kappa\theta}{(1-\theta\z)(1-\theta\w)}-\frac \theta{\z\w}.
$$
By \eqref{e2} we get, cf. \eqref{bess},
\begin{multline}\label{conf}
(i-j)\bK(i,j)\\=\frac{\theta}{(2\pi\sqrt{-1})^2}
\iint\limits_{|\z|=c_1<\frac 1{|\theta|},\,|\w|=c_2}
\Biggl(
\kappa\,\frac{(1-\theta\z)^{-\kappa-1}e^{-\theta\z^{-1}}}{\z^{i+1}}\,\frac
{(1-\theta\w)^{\kappa-1}
e^{\theta\w^{-1}}}{\w^{-j}}\\-
\frac{(1-\theta\z)^{-\kappa}e^{-\theta\z^{-1}}}{\z^{i+2}}\,
\frac{(1-\theta\w)^{\kappa}e^{\theta\w^{-1}}}{\w^{-j+1}}\Biggr)
{d\z}{d\w}.
\end{multline}
The contour integrals can be easily computed in terms of
the confluent hypergeometric function
$\Phi(a,c;x)={}_1F_1({a};{c};x)$: for any positive integer $m$
$$
\frac{1}{2\pi\sqrt{-1}}\int\limits_{|\z|=c_1<\frac 1{|\theta|}}
(1-\theta\z)^{-\alpha}e^{-\theta\z^{-1}}\frac{d\z}{\z^{m+1}}=
\frac{(\alpha)_m}{m!}\,e^{-\theta^2}\,\Phi(1-\alpha,m+1;\theta^2),
$$
$$
\frac{1}{2\pi\sqrt{-1}}\int\limits_{|\w|=c_2}
(1-\theta\w)^{\beta}e^{\theta\w^{-1}}\,\w^{m-1}{d\w}=
\frac{1}{m!}\,\Phi(-\beta,m+1;\theta^2)
$$
where $(a)_k=a(a+1)\cdots(a+k-1)$, $(a)_0=1$, is the Pochgammer symbol.

Hence, \eqref{conf} yields
\begin{equation}
\begin{gathered}\label{char}
\bK(i,j)=\frac{(\kappa)_{i+1}e^{-\theta^2}}{i-j}\,\Biggl(
\frac{\Phi(-\kappa,i+1;\theta^2)}{i!}\,
\frac{\Phi(-\kappa+1,j+2;\theta^2)}{(j+1)!}\\-
\frac{\Phi(-\kappa+1,i+2;\theta^2)}{(i+1)!}\,
\frac{\Phi(-\kappa,j+1;\theta^2)}{j!}\Biggr).
\end{gathered}
\end{equation}

As in the previous example, it can be shown that the diagonal entries
$\bK(i,i)$ can be obtained as the limit of the right--hand side
of \eqref{char} when $j\to i$.

The constant $Z$ in this example equals $\exp(\kappa\theta^2)$, and Theorem
\ref{main} gives
\begin{equation}\label{ex2}
D_n((1+\theta\z)^\kappa\exp(\theta \z^{-1}))=\exp(\kappa\theta^2)
\,\det (1-\bK)_{\ell_2(\{n,n+1,\dots\})}\,
\end{equation}
where $|\theta|<1$ and $\bK$ is as above.

Toeplitz determinants from the left--hand side of \eqref{ex2}
with positive integral $\kappa$ have been studied in \cite{TW}
in connection with asymptotics of the length of longest increasing
subsequences in random words.

On the other hand, if $\kappa$ is a positive integer,
the kernel $\bK$ described above
turns into a kernel conjugate to
the Christoffel--Darboux kernel of $\kappa$th order for
Charlier polynomials with parameter $\theta^2$, see \cite{HTF}
for definitions. This \emph{Charlier kernel} appeared in \cite{J}
in the study of the same combinatorial problem on random words.

For positive integral $\kappa$'s \eqref{ex2} follows from
the results of \cite{TW,J}.

The main result of the previous example (relation \eqref{ex1})
can be obtained from \eqref{ex2} by a simple limit transition.

Indeed, using the invariance of Toeplitz determinants with respect to
dilations, for any positive integer $k$ we can write that
\begin{equation}\label{lim1}
D_n((1+\widetilde\theta\z/k)^k\exp(\widetilde\theta\z^{-1}))=
D_n((1+\widetilde\theta\z/\sqrt{k})^k\exp(\widetilde\theta\z^{-1}/\sqrt{k}).
\end{equation}
As $k\to\infty$, the left--hand side of \eqref{lim1} tends to the left--hand
side of \eqref{ex1} with $\theta=\widetilde \theta$. On the other hand,
one can check that the Charlier kernel of $k$th order with parameter
$\theta^2=\widetilde \theta^2/{k}$ converges, as $k\to\infty$,
to the kernel of the previous example; \cite{J} contains a detailed
discussion of this limit transition, see also \cite{BO3}.
\end{Example}

\begin{Example}\label{example3} \rm 
Here we will be interested in Toeplitz determinants
with symbols of the form $\vph(\z)=(1+\xi_1\z)^z(1+\xi_2\z^{-1})^{z'}$
where $\xi_1,\,\xi_2,\,z,\,z'$ are complex parameters. Invariance of Toeplitz
determinants with respect to dilations implies that it suffices to consider
the case when $\xi_1=\xi_2=\xi$. Then
$$
V(\z)=z\ln(1+\xi\z)+z'\ln(1+\xi\z^{-1})
$$
and conditions of Theorem \ref{main} require that $|\xi|<1$.
We have
$$
\tV(\z)=-z\ln(1-\xi\z)+z'\ln(1-\xi\z^{-1}),\quad
\z\frac d{d\z}\tV(\z)=\frac{z\xi\z}{1-\xi\z}+\frac{z'\xi\z^{-1}}{1-\xi\z^{-1}},
$$
$$
\frac{\z\dfrac{d}{d\z} \tV(\z)- \w\dfrac{d}{d\w} \tV(\w)}{\z-\w}=
\frac{z\xi}{(1-\xi\z)(1-\xi\w)}-
\frac{z'\xi}{\z\w(1-\xi\z^{-1})(1-\xi\w^{-1})}.
$$
By \eqref{e2} we get, cf. \eqref{bess}, \eqref{conf},
\begin{multline}\label{gauss}
(i-j)\bK(i,j)=\frac 1{(2\pi\sqrt{-1})^2}\\ \times
\iint\limits_{|\xi|<|\z|=c_1,|\w|=c_2<\frac 1{|\xi|}}
\Biggl(
z\,
\frac{(1-\xi\z)^{-z-1}(1-\xi\z^{-1})^{z`}}{\z^{i+1}}\,
\frac{(1-\xi\w)^{z-1}(1-\xi\w^{-1})^{-z'}}{\w^{-j}}\\-
z'\,
\frac{(1-\xi\z)^{-z}(1-\xi\z^{-1})^{z'-1}}{\z^{i+2}}\,
\frac{(1-\xi\w)^{z}(1-\xi\w^{-1})^{-z'-1}}{\w^{-j+1}}\Biggr)
{d\z}{d\w}.
\end{multline}

This time the contour integrals can be expressed via the Gauss
hypergeometric function: for any integer $m\ge 0$
$$
\begin{aligned}
\frac 1{2\pi\sqrt{-1}}\int\limits_{|\xi|<|\z|=c_1<\frac 1{|\xi|}}
&(1-\xi\z)^{-\alpha}(1-\xi\z^{-1})^{\alpha'}\frac{d\z}{\z^{m+1}}\\&=
\frac{(\alpha)_m}{m!}\,\xi^m(1-\xi^2)^{-\alpha'}
F\left(1-\alpha,\alpha';m+1;\frac{\xi^2}{\xi^2-1}\right),
\end{aligned}
$$
$$
\begin{aligned}
\frac 1{2\pi\sqrt{-1}}\int\limits_{|\xi|<|\w|=c_2<\frac 1{|\xi|}}
&(1-\xi\w)^{\beta}(1-\xi\w^{-1})^{-\beta'}\,\w^{m-1}{d\w}\\&=
\frac{(\beta')_m}{m!}\,\xi^m(1-\xi^2)^{-\beta}
F\left(\beta,1-\beta';m+1;\frac{\xi^2}{\xi^2-1}\right).
\end{aligned}
$$

Substituting these formulas in \eqref{gauss}, we get
\begin{equation}\label{hyper}
\begin{gathered}
\bK(i,j)=\frac{(z)_{i+1}(z')_{j+1}}{i!j!}\,\xi^{i+j+2}(1-\xi^2)^{z+z'-1}
\\ \times \frac 1{i-j}\Biggl(
F\left(-z,-z';i+1;\frac{\xi^2}{\xi^2-1}\right)
\frac{F\left(1-z,1-z';j+2;\frac{\xi^2}{\xi^2-1}\right)}{j+1}\\-
\frac{F\left(1-z,1-z';i+2;\frac{\xi^2}{\xi^2-1}\right)}{i+2}
F\left(-z,-z';j+1;\frac{\xi^2}{\xi^2-1}\right).
\Biggr)
\end{gathered}
\end{equation}

As in two examples above, the diagonal values $\bK(i,i)$ can be computed
as limits of the right--hand side of \eqref{hyper} as $j\to i$.

The kernel thus obtained is a part of the \emph{hypergeometric kernel}
obtained in \cite{BO2}. The hypergeometric kernel describes
a remarkable 3--parametric
family of measures on partitions called \emph{z--measures} which are closely
related to so--called generalized regular
representations of the infinite symmetric group,
see \cite{KOV,BO1,BO2,BO3}. It lies on top of a hierarchy 
of kernels which are expressible in classical special functions,
see \cite{BO3} for details.

Let us compute the constant $Z$. We have
$$
v_k^+=\frac{z(-\xi)^k}{k},\quad
v_k^-=\frac{z'(-\xi)^k}{k}.
$$
Hence,
$$
Z=\exp\left(zz'\sum_{k=1}^\infty\frac{\xi^{2k}}{k}\right)=(1-\xi^2)^{-zz'}.
$$
Then Theorem \ref{main} gives
\begin{equation}\label{ex3}
D_n((1+\xi\z)^z(1+\xi\z^{-1})^{z'})=(1-\xi^2)^{-zz'}
\,\det (1-\bK)_{\ell_2(\{n,n+1,\dots\})}\,
\end{equation}
where $|\xi|<1$ and $\bK$ is as above.

When one of parameters $z,z'$, say $z$, is a positive integer,
the kernel $\bK$
turns into the Christoffel--Darboux kernel of order $z$ for Meixner polynomials
with parameters $(z'-z,\xi^2)$,
see \cite{BO2}, \S4, for details.

Formulas \eqref{ex1}, \eqref{ex2} can be obtained as certain
limits of \eqref{ex3}. To get \eqref{ex1} we just take
$\xi=\theta/k$, $z=z'=k$ with positive integers $k$ and
send $k$ to the infinity.

To get \eqref{ex2}
we employ the same trick as was
used in the previous example. Using invariance of Toeplitz determinants
with respect to dilations, we look at Toeplitz determinants
with symbols $(1+\theta\z)^{\kappa}(1+\theta\z^{-1}/k)^k$ for positive
integers $k$ and take the limit $k\to\infty$.

A detailed discussion of these limit transition in the language of
kernels and underlying combinatorial problems
can be found in \cite{J, BO3}.
\end{Example}

\noindent 
{\em E-mail:}\quad {\tt borodine@math.upenn.edu},
\quad {\tt okounkov@math.berkeley.edu}

\end{document}